\documentclass[a4paper]{article}
\usepackage{amsmath}
\usepackage{amssymb}
\usepackage{eucal}

\newtheorem{thm}{Theorem}[section]
\newtheorem{lem}[thm]{Lemma}
\newtheorem{cor}[thm]{Corollary}
\newtheorem{prop}[thm]{Proposition}
\newtheorem{rem}[thm]{Remark}
\newtheorem{defn}[thm]{Definition}
\newtheorem{eg}[thm]{Example}

\newenvironment{prf}{{\noindent Proof:\ }}{\hfill $\Box$\\ \smallskip}

\newcommand{\smnoind}{{\smallskip\noindent}}

\newcommand{\id}{{\rm id}}

\newcommand{\ti}{\tilde}

\newcommand{\cb}{{\rm CB}}

\newcommand{\abs}[1]{\left\vert#1\right\vert}
\newcommand{\K}{\mathcal{K}}
\newcommand{\BK}{\mathcal{L}_{\K(l^2)}}

\newcommand{\re}{{\rm reg}}

\newcommand{\uba}{\mathfrak{A}}

\newcommand{\e}{{\rm ess}}

\begin{document}
\title{Regular normed bimodules}
\author{Chi-Keung Ng\thanks{This work is supported by the National Natural Science Foundation of China (10371058).}}
\date{}
\maketitle
\begin{abstract}
In this article, we will give a characterization of Banach bimodules over $C^*$-algebras of compact operators that arises from operator spaces as well as a characterization of (F)-Banach bundles amongst all (H)-Banach bundles over a hyper-Stonian space. 
These two characterizations are concerned with whether certain natural map from a Banach bimodule to its canonical bidual is isometric (we call such bimodule \emph{regular}). 

\end{abstract}

\bigskip
\bigskip

\section*{Introduction}

\bigskip

The aim of this paper is to study duality theory for essential normed bimodules. 
Given a pre-$C^*$-algebra $A$ and an essential  
normed $A$-bimodule $X$, we would like to have a canonical definition of the dual object $X^s$ 
of $X$ which satisfies the following properties: 

\begin{enumerate}
\item $X^s$ is also an essential normed $A$-bimodule (i.e. the dual object is in the same category);

\item $X^s$ depends only on $X$ and $A$; 

\item when $A = \K(l^2)$ and $X$ is defined by an operator space, $X^s$ is defined by the corresponding dual operator space;

\item when $A$ is commutative, $X^s$ is the essential part of $\mathcal{L}_A(X,A)$ (i.e. the duality 
agree with the usual one for commutative algebras). 
\end{enumerate}

Let's forget about the norm structure for the moment and consider a  bimodule $M$ 
over a unital algebra $R$. 
The natural ``dual object'' $L_R(M,R)$ fails to be a $R$-bimodule if $R$ is not commutative. 
An easy way to rectify this situation is to ``add another copy of $R$'' and consider $L_R(M,\mathfrak{R})$ where $\mathfrak{R}$ is the algebraic tensor product $R\odot R$ together with the $R$-bimodule
structure: $a\cdot (b\otimes c) \cdot d = abd\otimes c$.
Therefore, $L_R(M,\mathfrak{R})$ becomes a $R$-bimodule (given by the bimodule structure on the second variable of $R\odot R$). 
However, when $R$ is commutative, 
$L_R(M,\underline{R}\odot R) \neq L_R(M,R)$ unless $R$ is the scalar field. 
A natural way to correct this is to replace $R\odot R$ with $R\odot_Z R$ (where $Z$ is the center of $R$). 

\medskip

We employ this simple idea in Section 1 to define the ``regular dual object'', $X^s$, of an essential 
normed $A$-bimodule $X$ (for technical reason, we will assume that $A$ has a contractive approximate identity and $A^2 = A$). 
There is a canonical contraction $\kappa_X: X \rightarrow X^{ss}$ (the dual of $X^s$). 
In general, $\kappa_X$ is not an isometry and $X$ is called \emph{regular}
if $\kappa_X$ happens to be an isometry. 
It is easy to see that $X^s$ is always regular and so, $\kappa_X(X)$ is called the 
\emph{regularization} of $X$. 
Regular bimodule are thought to be nice object because of the results in Sections 2 and 3. 
It is natural to ask whether one can give a canonical characterisation of regularity without 
explicitly involving the duality. 
We will show that every regular bimodule will satisfy certain properties (known as \emph{pseudo-regularity}) which do not involve the dual object explicitly. 

\medskip

In section 2, we will consider the situation when $A = \bigoplus_{\lambda\in \Lambda} \K(H_\lambda)$ where 
$\Lambda$ is any index set and $H_\lambda$ is a Hilbert space. 
We will show that regularity and pseudo-regularity coincide in this case and will characterize regular Banach  
$A$-bimodules. 

\medskip

In section 3, we will consider the situation when $A$ is a commutative von Neumann algebra. 
In this case, pseudo-regular bimodules correspond to 
(H)-Banach bundles and we will show that regular bimodules correspond to  
(F)-Banach bundles. 
Hence, the regularization process gives us a canonical way to obtain an (F)-Banach bundle from any 
given Banach bundle (in particular, from an (H)-Banach bundle). 

\noindent \emph{Acknowledgment:} 
The author would like to thank Prof. Effros and Prof. Ruan for some helpful discussions on an earlier version of this article (namely, \cite{Ng-OS}). 
The author would also like to thank the referee for helpful comments that lead to simplifications of the arguments of some results and for informing us about \cite{Mag-dual}. 

\bigskip
\bigskip

\section{Duality and Regularity of Banach bimodules}

\bigskip

Throughout this section, $A$ is a pre-$C^*$-algebra containing a contractive  
approximate identity $\{f_i\}_{i\in I}$ such that 
$A\cdot A = A$ (i.e. any element of $A$ is a finite sum of elements of the form $ab$ where $a,b\in A$). 
We denote by $\bar A$ the completion of $A$ and recall that an $A$-bimodule $X$ is \emph{(algebraically) essential} if $X = A\cdot X \cdot A$ 
(where $A\cdot X \cdot A$ is the linear span of elements of the form $a\cdot x\cdot b$ with $a,b\in A$ and $x\in X$). 
For any $A$-module $X$, we denote by $X_E$ the \emph{essential part} $A\cdot X \cdot A$ of $X$. 

\medskip

Let $B$ be a $C^*$-subalgebra of $M(\bar A)$ (the multiplier algebra of $\bar A$) and 
$Z_A \cong C(\Omega)$ be the 
center of $M(\bar A)$ (where $\Omega$ is a compact Hausdorff space). 
By \cite{Bl}, there exists a $C^*$-semi-norm $\|\cdot \|_m$ on the algebraic $Z_A$-tensor 
product $\bar A\odot_{Z_A} B$ which is minimum in some sense (see \cite[2.8]{Bl}). 
As in \cite{Bl}, we denote by $\bar A \stackrel{m}{\otimes}_{Z_A} B$ the Hausdorff 
completion of $\bar A\odot_{Z_A} B$ under $\| \cdot\|_m$. 

\medskip

From now on, we will denote by $a\circledast_A 1$ and $a\circledast_A b$ 
(or simply $a\circledast 1$ and $a\circledast b$) the canonical images of $a\in A$ and $a\otimes_{Z_A} b\in (A\odot_{Z_A} B)$ 
in $M(\bar A \stackrel{m}{\otimes}_{Z_A} B)$ and $\bar A \stackrel{m}{\otimes}_{Z_A} B$ 
respectively.
Note that the map that sends $a\in A$ to 
$a\circledast 1\in M(\bar A \stackrel{m}{\otimes}_{Z_A} B)$ is a $*$-homomorphism.

\medskip

\begin{lem}
\label{tp-o-Z}
\ \ $\bar A \stackrel{m}{\otimes}_{Z_A} B$\ \  is a normed $A$-bimodule under the multiplication 
$c\cdot (a\circledast b)\cdot d = cad\circledast b$ 
and $\{ f_i\}_{i\in I}$ is an approximate identity in $A$ for the $A$-bimodule $\bar A \stackrel{m}{\otimes}_{Z_A} B$. 
Similarly, if $\{g_j\}$ is a bounded approximate identity of $B$, then both 
$(1\circledast g_j)\alpha$ and $\alpha(1\circledast g_j)$ converge to $\alpha$ for any 
$\alpha\in \bar A \stackrel{m}{\otimes}_{Z_A} B$. 
\end{lem}

\medskip

The above lemma may be used implicitly throughout this article. 
In the following we denote by $\uba$
the normed $A$-bimodule $\bar A \stackrel{m}{\otimes}_{Z_A} \bar A$ 
(the $A$-bimodule structure is as given in Lemma \ref{tp-o-Z}). 

\medskip

We will now construct the dual bimodule of an essential normed  $A$-bimodule. 
We have already stated in the introduction what are expected for dual bimodules. 
For any essential normed $A$-bimodule $X$, we denote by $\mathcal{L}_A(X, \uba)$ the space of all continuous $A$-bimodule maps from $X$ to $\uba$.
Note that 
$\mathcal{L}_A(X, \uba) = \mathcal{L}_A(X, A\cdot \uba \cdot A)$ is a $A$-bimodule with multiplications given by 
$(a\cdot T\cdot b)(x) = (1\circledast a)T(x)(1\circledast b)$ ($T\in \mathcal{L}_A(X, \uba)$; $a,b\in A$; $x\in X$). 
We set  
$$X^s \ :=\ \mathcal{L}_A(X, \uba)_E$$
(the essential part of $\mathcal{L}_A(X, \uba)$) and call it the \emph{regular dual} of $X$. 

\medskip

\medskip

It is clear that $X^s$ is an essential normed $A$-bimodule 
(with the canonical norm on $\mathcal{L}_A(X, \uba)$) and there exists a contraction 
\ $\kappa_X : X\rightarrow X^{ss}$\ given by 
$\kappa_X(x)(\varphi) = \varphi(x)^{(12)}$ 
(where $(12)$ is the flip of the two variables in $\uba$ -- it is not hard to see from the definition of $\stackrel{m}{\otimes}_{Z_A}$ that such 
flip map exists). 

\medskip

\begin{defn}
\label{def-reg}
Let $X$ be an essential normed $A$-bimodule. 
An element $x\in X$ is said to be \emph{regular in $X$} if for any $\epsilon > 0$, there 
exists $T\in X^s$ such that $\|T\| \leq 1$ and $\|x\| \leq \|T(x)\| + \epsilon$. 
We say that $X$ is \emph{regular} if $\kappa_X$ is an isometry.
Moreover, the closure of $\kappa_X(X)$ in $X^{ss}$ together with the induced norm is called the \emph{regularization} of $X$ and is denoted by $X_\re$.
\end{defn}

\medskip

In \cite{Pop}, a type of dual bimodule, $X^\dagger$, was introduced. 
However, $X^\dagger$ is in 
general not a Banach $A$-bimodule. 
Furthermore, the regularity defined 
using $X^\dagger$ (i.e. $X\rightarrow (X^\dagger)^\dagger$ being isometric) is in general strictly weaker than ours. 
In particular, it will not give the relation between (H)-Banach bundles and 
(F)-Banach bundles as obtained using our notion of regularization (see Section 3 below). 

\medskip 
In \cite{Na}, yet another type of dual bimodule $X^*_\mathcal{D}$ was introduced. 
If $A \subseteq \mathcal{L}(H)$ such that $A''$ is standard in $\mathcal{L}(H)$, then $X^*_{\mathcal{L}(H)}$ is very similar to $X^\dagger$ (except that $\rm Hom$ in \cite{Pop} are bounded maps while $\rm Hom$ in \cite{Na} are completely bounded). 
If $A = \K(H)$ or if $A$ is commutative, then $X^*_\uba = X^s$ (in both cases, elements in $X^s$ are automatically completely bounded) but it seems unlikely that one can use any result in \cite{Na} to shorten the proofs in this paper. 

\medskip

\begin{rem}
\label{complete}
\ Let $X$ be an essential normed $A$-bimodule and $\bar X$ be its  
completion. 

\smnoind
(a)\ If $A$ is a $C^*$-algebra, then $X^s$ is closed in $\mathcal{L}_A(X, \uba)$ 
(because of the the Cohen factorization theorem). 

\smnoind
(b)\ $\bar X$ is an essential Banach $\bar A$-bimodule and $X^s$ is dense in $\bar X^s$ (since 
$\mathcal{L}_A(X, \uba) = \mathcal{L}_{\bar A}(\bar X, \uba)$).
Thus, $X$ is regular if and only if $\bar X$ is regular. 

\smnoind
(c)\ $X$ is regular if and only if for any 
$x\in X$ and any $\epsilon >0$, there exists 
$T\in \mathcal{L}_A(X,\uba) = \mathcal{L}_A(X, A\cdot \uba \cdot A)$ such that 
\ $\| T\| \leq 1$\ and \ $\|x\| < \|T(x)\| +\epsilon$ (note that $\|f_i\cdot T(x) \cdot f_i - T(x)\| < \epsilon /2$ for large enough $i$ because of Lemma \ref{tp-o-Z}). 
Consequently, essential submodule of a regular bimodule is again regular. 

\smnoind
(d) Since $(\kappa_X)^s\circ \kappa_{X^s} = \id_{X^s}$ and both $\kappa_{X^s}$ and $(\kappa_X)^s: (X^{ss})^s \rightarrow X^s$ are contractions, $\kappa_{X^s}$ is always an isometry. 
Therefore, $X^s$ is regular and so is $X_{\rm reg}$ (by 
part (b)). 
Moreover, if $Y$ is an regular normed $A$-bimodule and 
$T\in \mathcal{L}_A(X,Y)$, there exists $T_\re\in \mathcal{L}_{\bar A}(X_\re,Y)$ such that 
$T = T_\re\circ \kappa_X$. 

\smnoind
\end{rem}

\medskip

It is natural to ask if one can characterise regularity without finding the  
regular dual. 
In some cases, this can be done using the notion of pseudo- 
regularity as defined in the following (although, this is not the case in general; see e.g. Section 3). 

\medskip

\begin{defn}
\label{abs-conv}
(a) A semi-norm $p$ on a $A$-bimodule $X$ is said to be  
\emph{absolutely $A$-convex} if for any $a_1,...,a_n,b_1,...,b_n\in A$ and any $x_1,...,x_n\in X$, 
$$p\left(\sum_{i=1}^n a_i x_i b_i\right) \ \leq \  \sqrt{\left\|\sum_{i=1}^n a_ia_i^* \right\|} \max_{i=1,...,n} p(x_i) \sqrt{\left\| \sum_{i=1}^n b_i^*b_i \right\|}.$$ 
Moreover, an essential normed $A$-bimodule is said to be \emph{absolutely $A$-convex} if its norm is absolutely $A$-convex. 

\smnoind
(b) An essential normed $A$-bimodule $X$ is said to be \emph{pseudo-regular} if it is  
absolutely 
$A$-convex and its completion $\bar X$ (which is a unital $M(\bar A)$-bimodule in the canonical way) is a commutative $Z_A$-bimodule. 
\end{defn}

\medskip

If $A$ is a $C^*$-algebra and $X$ is an essential Banach $A$-bimodule, then  
\cite[Proposition 2.2]{Pop} tells us that one only needs to consider $n=2$ in the  
definition of absolute $A$-convexity

\medskip

\begin{eg}
\label{reg=conv}
(a) If $A = \bigoplus_{i\in I}^{c_0} \K(H_i)$ ($c_0$-direct sum), then we have  $\uba = \bigoplus_{i\in I}^{c_0} \K(H_i)\otimes \K(H_i)$.

\smnoind
(b) If $A = M_\infty$ (the space of all infinite matrices with finite numbers of non-zero entries, considered as a subspace of $\K(l^2)$), then $\uba = \K(l^2)\otimes \K(l^2)$. 

\smnoind
(c) If $A = C_0(\Omega)$ for some locally compact space $\Omega$, then 
$\uba = C_0(\Omega)$. 

\smnoind
(d) Let $H$ be an infinite dimensional Hilbert space, $\mathbf{W}$ be an operator space, 
$X := \mathbf{W}\check\otimes \K(H)$ (spatial tensor product) and $X^\# := \mathcal{L}_{\K(H)}(X; \K(H)\check\otimes \K(H))$.
It is clear that for any $T\in \cb(\mathbf{W}; \K(H))$, we have $T\otimes \id_{\K(H)} \in X^\#$. 
Conversely, any $\varphi \in X^\#$ restricts to a map $$\varphi_0\ \in\ \mathcal{L}_{\K(l^2)}(\mathbf{W}\check\otimes \K(l^2); \K(H)\check\otimes \K(l^2))\ \cong\  \cb(\mathbf{W};\K(H))$$ (see e.g. \cite[1.2]{Ng-cohom}).
If $T\in \cb(\mathbf{W};\K(H))$ is the corresponding element of $\varphi_0$, then clearly $\varphi = T\otimes \id_{\K(H)}$. 
Consequently, $X^s \cong \cb(\mathbf{W}; \K(H))_E$ (see the first 
paragraph). 

\smnoind
(e) For any Hilbert space $H$, one can consider $\K(H)^*$ as an essential normed $\K(H)$-bimodule. 
In this case, $(\K(H)^*)_\re = \K(H)$. 
In fact, $\K(H)^*\rightarrow (\K(H)^*)_\re$ is 
the identification of $\K(H)^*$ as the set of trace-class operators. 

\smnoind
(f) Suppose that $A = \K(l^2)$ or $M_\infty$ and $X$ is an essential normed $A$-bimodule.
Then a closed subset $D\subseteq X$ is absolutely $A$-convex if and only if 
for any disjoint projections $p,q\in A$ and any $a \in A$ with $\|a\|\leq 1$, 
we have $p\cdot D\cdot p + q\cdot D\cdot q \subseteq D$, $a\cdot D\subseteq D$ and 
$D\cdot a \subseteq D$.
In fact, the case of $A=M_\infty$ is more or less the same as \cite[3.2]{EW} and the case of $A=\K(l^2)$ follows from some completion arguments. 
\end{eg}

\medskip

\begin{prop}
\label{ub>abconv}
If $X$ is a regular normed $A$-bimodule, $X$ is pseudo-regular. 
\end{prop}
\noindent
\begin{prf}
\ Let $U$ and $V$ be the closed unit balls of $X$ and $X^s$ respectively. 
The  
regularity of $X$ means that $U = \{ x\in X: \|\varphi(x)\|\leq 1$ for any $\varphi\in V\}$. 
Therefore, $X$ is absolutely $A$-convex (note that the norm on $\uba$ is absolutely $A$-convex). 
For any $\varphi\in \bar X^s$, $z\in Z_A$ and $x = \sum_{k=1}^n a_kx_kb_k \in X$ ($a_k, b_k\in A$ and 
$x_k\in X$), we have 
$$\varphi(z\cdot x)\ =\ \sum_{k=1}^n (za_k\circledast 1)\varphi(x_k)(b_k\circledast 1) 
\ =\ \sum_{k=1}^n (a_k\circledast 1)\varphi(x_k)(b_kz\circledast 1) 
\ =\ \varphi(x\cdot z)$$ 
and so $z\cdot x = x\cdot z$ (as $\bar X^s$ separates points of $\bar X$ by Remark \ref{complete}(c)). 
Thus, $X$ is pseudo-regular because the multiplications are continuous. 
\end{prf}

\bigskip
\bigskip

\section{The case when $A = \bigoplus_{i\in I}\mathcal{K}(H_i)$}

\bigskip

We will first consider the case when $A = \mathcal{K}(H)$ (where $H$ is a Hilbert space). 
In the following, $\check \otimes$ is the spatial tensor product of two operator spaces.

\medskip

\pagebreak

\begin{thm}
\label{os}
Let $X$ be an essential Banach $\K(H)$-bimodule.
The following statements are equivalent.

\smnoind
(1) $X$ is regular.

\smnoind
(2) $X$ is pseudo-regular. 

\smnoind
(3) There exists a complete operator space $\mathbf{W}$ such that $X = \mathbf{W} \check \otimes \K(H)$.
\end{thm}
\noindent
\begin{prf}
(1)$\Rightarrow$(2). 
This follows from Proposition \ref{ub>abconv}(b). 

\smnoind
(2)$\Rightarrow$(3). 
By the theorem in p.333 of \cite{Pop}, there exists a Hilbert space $K$, a non-degenerate $*$-representation $\pi$ of $\K(H)$ on $K$ as well as an isometry $J: X \rightarrow \mathcal{L}(K)$ such that $J(axb) = \pi(a)J(x)\pi(b)$ ($a,b\in \K(H);x\in X$). 
Note that there is a Hilbert space $E$ such that $K \cong E \otimes H$ as well as $\pi(a) = 1\otimes a$ and so we can assume $\pi = 1\otimes \id$ and $J = \id$. 
Let $\mathbf{W}:=\{y\in \mathcal{L}(E): y\otimes a \in X {\rm\ for\ any\ } a\in \K(H)\}$. 
Clearly, $\mathbf{W}\check \otimes \K(H) \subseteq X$. 
Suppose $\{\xi_i\}_{i\in I}$ is an orthonormal basis for $H$. 
For any $i,j\in I$, we set  
$$\theta_{i,j}(\zeta)\ :=\ \xi_i \langle \xi_j, \zeta \rangle 
\quad {\rm and} \quad \omega_{i,j}(t)\ :=\ \langle \xi_i, t \xi_j\rangle \qquad (\zeta \in H; t\in \mathcal{L}(H)).$$ 
Then $(\id\otimes \omega_{j,k})(x)\otimes \theta_{i,l} = (1\otimes \theta_{i,j})x(1\otimes \theta_{k,l}) \in X$ ($x\in X;i,j,k,l\in I$) and so $(\id\otimes \omega_{j,k})(x)\in \mathbf{W}$. 
Furthermore, since $X$ is essential, certain finite sums of elements of the form 
$$(1\otimes \theta_{i,i})x(1\otimes \theta_{j,j})\ =\ (\id\otimes \omega_{i,j})(x)\otimes \theta_{i,j}\ \in \mathbf{W} \check \otimes \K(H)$$ 
converge to $x$ in norm and so $X = \mathbf{W} \check \otimes \K(H)$ as required. 

\smnoind
(3)$\Rightarrow$(1). 
If $\dim H = n$, then $X \cong M_n(W)$ as normed $M_n$-bimodules and this implication follows from \cite[2.3.4]{ER} (note that $\cb(\mathbf{W};M_n) \cong X^s$). 
On the other hand, if $H$ is infinite dimensional, then $(\mathbf{V}\check\otimes \K(H))^s \cong \cb(\mathbf{V}; \K(H))_{E}$ for any operator space $\mathbf{V}$ (by Example \ref{reg=conv}(d)). 
Therefore, we have:
\begin{eqnarray*}
\mathbf{W}\check\otimes \K(H) \ \subseteq \ \cb(\mathbf{W}^*;\K(H))_E 
& = & (\mathbf{W}^*\check \otimes \K(H))^s \\ 
& = & (\cb(\mathbf{W}; \K(H))_E)^s \ = \ (\mathbf{W}\check\otimes \K(H))^{ss}  
\end{eqnarray*}
($\mathbf{W}^*\check \otimes \K(H) = \cb(\mathbf{W}; \K(H))_E$ because there is an approximate unit in $\K(H)$  
consisting of finite rank projections). 
It is not hard to check that the above  
embedding is precisely $\kappa_X$ and thus $X$ is regular. 
\end{prf}

\medskip

\begin{rem}
\label{incomplete}
(a) The equivalence of (2) and (3) is probably known (e.g. one can use \cite[2.1]{Mag} and \cite[Proposition 3.3]{Web} to obtain this in the case of $H = l^2$). 
However, we decided to give a proof here for clarity and completeness. 

\smnoind
(b) Let $\{\xi_i\}_{i\in I}$ be an othonormal basis for $H$ and $A$ be the linear span of  
$\{ \theta_{i,j}: i,j\in I\}$. 
One can use the completion consideration in Remark \ref{complete} to obtain a similar result as the above theorem for $A$. 
In fact, there is also an elementary proof for this fact (without using the theorem in \cite{Pop}) but such a proof is much more lengthy. 
\end{rem}

Suppose that $A$ is the $c_0$-direct sum $\bigoplus_{\lambda\in \Lambda}^{c_0} \mathcal{K}(H_\lambda)$ and $d_\lambda\in A$ corresponds to the identity in $\mathcal{L}(H_\lambda)$. 
Then $Z_A = c_0(\Lambda)$ and $\uba = \bigoplus_{\lambda\in \Lambda}^{c_0} \mathcal{K}(H_\lambda)\otimes \mathcal{K}(H_\lambda)$. 
Let $X$ be a pseudo-regular Banach $A$-bimodule. 
Then it is not hard to see that $X_\lambda := d_\lambda\cdot X$ is a regular Banach $\K(H_\lambda)$-bimodule and 
$X$ is the $c_0$-directed sum $\bigoplus_{\lambda\in \Lambda}^{c_0} X_\lambda$. 
Using this, one can check easily that $X$ is also regular. 
Thus, we have the following theorem. 

\medskip

\begin{thm}
\label{pse-reg=reg}
Let $\Lambda$ be an index set and $H_\lambda$ be a Hilbert space for any  
$\lambda\in \Lambda$. 
Suppose that $A = \bigoplus_{\lambda\in \Lambda}^{c_0} \mathcal{K}(H_\lambda)$ and $X$ is a pseudo-regular Banach  
$A$-bimodule. 
Then $X$ is regular and there exists a family of operator spaces $\{\mathbf{W}_\lambda\}_{\lambda\in \Lambda}$ such that $X = \bigoplus_{\lambda\in \Lambda}^{c_0} \mathbf{W}_\lambda\check \otimes \K(H_\lambda)$. 
\end{thm}

\medskip

The first two parts of the following corollary follow easily from Theorem \ref{os} (or more precisely, Remark \ref{incomplete}(b)) and the final part follows from the above theorem. 

\medskip

\begin{cor}
Let $X$ be an essential normed $A$-bimodule.

\smnoind
(a) Suppose that $A = \bigcup_{n\in \mathbb{N}} M_n$. 
Then $X$ is regular if and only if there exists an operator space $\mathbf{W}$ such that $X = \bigcup_{n\in \mathbb{N}} M_n(W)$. 

\smnoind
(b) Suppose that $A = M_n$. 
Then $X$ is regular if and only if there exists an operator space $\mathbf{Y}$ such that $X \cong M_n(Y)$ under the norm induced by the operator space structure of $\mathbf{Y}$. 

\smnoind
(c) Suppose that $A = c_0$. 
Then $X$ is regular if and only if there exists a sequence of Banach spaces $\{X_k\}$ such that $X$ is a normed $c_0$-submodule of $\bigoplus_{\lambda\in \Lambda}^{c_0} X_k$. 
\smnoind
\end{cor}

\medskip

In the remainder of this section, we will give two remarks concerning the case when $A = \K(l^2)$. 
First of all, Theorem \ref{os} allows us to detect some  
hidden operator space structures. 
For example, if $\mathbf{V}$ is a complete operator space, then 
any essential Banach $\K(l^2)$-submodule of $\K(l^2)\check\otimes \mathbf{V}$ is of the form
$\K(l^2)\check\otimes \mathbf{U}$ for some operator subspace $\mathbf{U}$ of $\mathbf{V}$. 
As for another example, if $Y$ is an essential  
operator $A$-bimodule of a $C^*$-algebra $A$, one can use the $\K(l^2)$-bimodule  
approach to show the existence of a canonical operator space structure on the space of double centralizers $M_A(Y)$ that turns it into a unital operator  
$M(A)$-bimodule in a canonical way (see e.g. \cite[p.310]{Ng-comod}).

\medskip

Secondly, ``regularization'' is a process that produces a canonical complete operator space from any essential Banach $\K(l^2)$-bimodule. 
The following  
corollary shows that it is actually a left adjoint of the forgetful 
functor from the category of complete operator spaces to the category of essential Banach  
$\K(l^2)$-bimodules 
(note that if $\mathbf{W}$ is the operator space such that $X_\re = \mathbf{W}\check\otimes \K(l^2)$, then the following corollary shows that $\cb(\mathbf{W}, \mathbf{V}) \cong \BK(X, \K(l^2)\check\otimes \mathbf{V})$  
canonically). 

\medskip

\begin{cor}
\label{reg}
\ Let $X$ and $Y$ be essential Banach $\K(l^2)$-bimodules. 
Any  
$\varphi\in\BK(X,Y)$ induces a map $\varphi_\re\in\BK(X_\re,Y_\re)$ such that $$\varphi_\re \circ \kappa_X = \kappa_Y\circ \varphi \quad {\rm and} \quad \|\varphi_\re\| \leq \| \varphi\|.$$ 
If, in addition, $Y$ is regular, then $\| \varphi_\re\| = \|\varphi\|$. 
Consequently, the canonical map, $\hat \kappa_X: \BK(X_\re, \K(l^2)\check\otimes \mathbf{V}) \rightarrow \BK(X, \K(l^2)\check\otimes \mathbf{V})$ is an isometry for any complete operator space $\mathbf{V}$. 
\end{cor}
\noindent
\begin{prf}
Consider $\varphi^s:Y^s \rightarrow X^s$ given by $\varphi^s(f) = f\circ \varphi$. 
It is easy to see that $\|\varphi^s\| \leq \| \varphi\|$. 
Hence we have a bounded $\K(l^2)$-bimodule map $\varphi^{ss}: X^{ss}\rightarrow Y^{ss}$ such that 
$$\varphi^{ss}\circ \kappa_X = \kappa_Y \circ \varphi$$ 
and $\|\varphi^{ss}\|\leq \|\varphi^s\| \leq \|\varphi\|$. 
Now, the restriction of $\varphi^{ss}$ on $X_\re$ is the required map $\varphi_\re$. 
Finally, if $Y$ is regular, $\kappa_Y$ is an isometry and so 
$\|\varphi\|\leq \|\varphi_\re\|\|\kappa_X\| \leq \|\varphi_\re\|$. 
\end{prf}

\bigskip
\bigskip

\section{The case when $A$ is a commutative von Neumann algebras}

\bigskip

Throughout this section, $\Omega$ is a compact Hausdorff space and $X$ is an  
essential Banach $C(\Omega)$-module (i.e. commutative Banach $C(\Omega)$-bimodule). 
For any $x\in X$, we denote by $X(x)$ the closed $C(\Omega)$-submodule $\overline{C(\Omega)\cdot x}$.

\medskip

As noted in \cite{Pop}, $X$ is pseudo-regular if and only if it is a $C(\Omega)$-convex  
module in the sense of \cite[p.40]{DG}. 
Therefore, it is the case if and only if $X$ is the space of continuous sections of an 
(H)-Banach bundle (see \cite[p.8]{DG} and \cite[2.5]{DG}).
Let us first give the following (probably well known) lemma. 

\medskip

\begin{lem}
\label{cont=leq}
Let $f\in C(\Omega)_+$ and $h: \Omega \rightarrow \mathbb{R}_+$ be an upper semi-continuous function. 
Then $f \leq h$ if and only if $\|gf\| \leq \|gh\| := \sup_{\omega\in \Omega} g(\omega)h(\omega)$ for any $g\in C(\Omega)_+$. 
\end{lem}
\begin{prf}
The necessity is clear. 
Suppose that there exist $\omega_0 \in \Omega$ and $r\in \mathbb{R}_+$ such that
$f(\omega_0) > r > h(\omega_0)$. 
Then $W = \{ \omega\in \Omega: h(\omega) < r < f(\omega)\}$ 
is an open set containing $\omega_0$. 
If $g$ is a continuous function from $\Omega$ to $[0,1]$ such that $0\leq g\leq 1$, $g(\omega_0) = 1$ and $g$ vanishes
outside $W$, then $\|gf\| > r > \|gh\|$. 
\end{prf}

\medskip

\begin{rem}
For any function $h: \Omega \rightarrow \mathbb{R}_+$, we define $\|h\| := \sup_{\omega\in \Omega} h(\omega)$ and 
$$\|h\|_{e} \ := \ \inf \left\{ \sup_{\omega\in \Delta} h(\omega): \Delta {\rm \ is\ an\ open\ 
dense\ subset\ of\ }\Omega \right\}.$$ 
If $h$ is upper-semi-continuous, then 
$$\|h\|_{e} = \inf \left\{ \sup_{\omega\in \Xi} h(\omega): \Xi {\rm \ is\ a\ dense\ subset\ of\ } \Omega \right\}.$$
\end{rem}

\medskip

\begin{prop}
\label{equiv-loc-reg}
Suppose that $X$ is a $C(\Omega)$-convex Banach module and $x\in X$. 
Define $\|x\|_e := \| \abs{x}\|_{e}$ (where $\abs{x}(\omega)=\|x(\omega)\|$). 
Then the following statements are equivalent.

\smnoind
(i) $\| x \| = \| x\|_{e}$. 

\smnoind
(ii) For any $\epsilon > 0$, there exists $f\in C(\Omega)_+$ such that $f\leq \abs{x}$ and 
$\| x \| \leq \|f\| + \epsilon$. 

\smnoind
(iii) $x$ is regular in $X(x) = \overline{C(\Omega)\cdot x}$ (see Definition \ref{def-reg}).

\smnoind
Consequently, if $X$ is regular, then $\|x\| = \|x\|_{e}$ for any $x\in X$. 
\end{prop}
\begin{prf}
(i)$\Rightarrow$(ii) 
Since $G := \{\omega \in \Omega: \abs{x}(\omega) \geq \|x\|_e - \epsilon\}$ is a closed set in $\Omega$ (as $\abs{x}$ is upper semi-continuous), $G$ contains an open set $V$ (otherwise, the open set 
$\{\omega\in \Omega: \abs{x}(\omega) < \|x\|_e - \epsilon \}$
is dense which contradicts the definition of $\|x\|_e$). 
Take any $\omega_0\in V$. 
Let $f\in C(\Omega)$ be such that $0\leq f(\omega) \leq \|x\|_e - \epsilon$ ($\omega\in \Omega$), 
$f(\omega_0) = \|x\|_e - \epsilon$ and $f$ vanishes outside $V$. 
Then clearly $f\leq \abs{x}$ and $\|x\|_e = \|f\| + \epsilon$. 

\smnoind
(ii)$\Rightarrow$(iii) 
For any $\epsilon > 0$, let $f$ be the function as given in statement (ii). 
We first show that $\varphi: X(x) \rightarrow C(\Omega)$ given by $\varphi(g\cdot x) = gf$ is well defined. 
Suppose that $g\in C(\Omega)$ such that $g(\omega) x(\omega) = 0$ for any $\omega\in \Omega$. 
If $g(\omega) \neq 0$, then $x(\omega) = 0$ and so $f(\omega) = 0$ which implies that $gf = 0$. 
Thus, $\varphi\in X(x)^s$ is a well defined contraction such that $\|x\| \leq \|\varphi(x)\| + \epsilon$. 

\smnoind
(iii)$\Rightarrow$(i) 
It is clear that $\|x\|_{e} \leq \|x\|$. 
For any $\epsilon > 0$,  let $\varphi\in X(x)^s$ such that $\|\varphi\| \leq 1$ and $\|x\| \leq 
\|\varphi(x)\| + \epsilon$. 
Put $f = \abs{\varphi(x)}\in C(\Omega)_+$.
Then for any $g\in C(\Omega)_+$, we have 
$(g\abs{\varphi(x)})(\omega) = \abs{\varphi(g\cdot x)}(\omega)$ for any $\omega\in \Omega$ and so, 
\begin{eqnarray*}
\|gf\|\ =\ \sup_{\omega \in \Omega} \abs{\varphi(g\cdot x)(\omega)}
& = & \|\varphi(g\cdot x)\|\\ 
& \leq & \|g\cdot x\|\ =\ \sup_{\omega \in \Omega} g(\omega) \|x(\omega)\|\ =\ \|g\abs{x}\|.
\end{eqnarray*}
Hence $f \leq \abs{x}$ (by Lemma \ref{cont=leq}) and $\| f\| = \| f\|_{e} \leq \| x \|_e$. 
Thus, 
$$\| x \|\ \leq\ \|\varphi(x)\| + \epsilon\ =\ \|f\| + \epsilon \leq \| x \|_{e} + \epsilon.$$
\end{prf}

\medskip

\begin{lem}
\label{ess-semi-norm}
Let $X$ be a $C(\Omega)$-convex Banach module. 
The map $x \mapsto \|x\|_e$ is an absolutely $C(\Omega)$-convex seminorm on $X$.
\end{lem}
\begin{prf}
Let $f_1,f_2\in C(\Omega)_+$ with $\|f_1 + f_2\| \leq 1$ and $x_1,x_2\in X$ with 
$\|x_1\|_e, \|x_2\|_e \leq 1$. 
For any $\epsilon > 0$, there exist open dense subsets $\Delta_1$ and $\Delta_2$ such that $\sup_{\omega\in \Delta_i} \|x_i(\omega)\| < 1 + \epsilon$ ($i=1,2$). 
If $\Delta = \Delta_1 \cap \Delta_2$, 
\begin{eqnarray*}
\|f_1\cdot x_1 + f_2\cdot x_2\|_e
& \leq & \sup_{\omega\in \Delta} \|f_1(\omega)x_1(\omega) + f_2(\omega)x_2(\omega)\| \\
& \leq & \sup_{\omega\in \Delta} f_1(\omega)\|x_1(\omega)\| + f_2(\omega)\|x_2(\omega)\|
\ = \ 1 + \epsilon.
\end{eqnarray*}
\end{prf}

\begin{rem}
\label{t2>norm}
$\|\cdot\|_e$ is a norm if the underlying topology of the (H)-Banach bundle $(p,E,\Omega)$ associated 
with $X$ is Hausdorff. 
In fact, suppose that $y\in X$ such that $\|y\|_e =0$. 
For any $n\in \mathbb{N}$, the open set $\{ \omega: \|y(\omega)\| < 1/n\}$ is dense in $\Omega$. 
Therefore, by the Baire's Category theorem, 
$$K_y \ := \ \{\omega\in \Omega: y(\omega) = 0_\omega\} 
\ = \ \bigcap_{n\in \mathbb{N}} \{ \omega: \|y(\omega)\| < 1/n\}$$ 
is dense in $\Omega$ (where $0_\omega$ is the zero of the fibre at $\omega$). 
Consider the map $j: \Omega \rightarrow E$ defined by $j(\omega) = 0_\omega$. 
By condition (4) of \cite[1.1]{DG}, we see that $j$ is a continuous map and so $j(\Omega)$ is compact in the Hausdorff space $E$. 
Thus, $K_y = y^{-1}(j(\Omega))$ is also closed in $\Omega$ and hence $K_y = \Omega$. 
This shows that $y \equiv 0$. 
Thus, $\|\cdot\|_e$ is a norm on $X$. 
\end{rem}

\medskip

In the following, we denote by $X_\e$ the completion of $(X/N, \|\cdot\|_e)$ (where $N = \{x\in X: \|x\|_e = 0 \}$).  
A natural question is whether $X = X_\e$ for any 
absolutely $C(\Omega)$-convex Banach module $X$. 
This is of course, true if $\Omega$ is a 
finite set. 
The following example shows that it is not the case in general. 

\medskip

\begin{eg}
\label{eg-h}
Let $\Omega$ be a compact Hausdorff space with a non-isolated point $\omega\in \Omega$. 
One can turn $\mathbb{C}$ into a Banach $C(\Omega)$-module, denoted by $X$, through the multiplication $f\cdot r = f(\omega)r$. 
It is not hard to check that $X$ is $C(\Omega)$-convex and so $X = \Gamma(E)$ for an (H)-Banach bundle $E$ over $\Omega$. 
By the construction in 
\cite[p.35-36]{DG}, we see that the fibre $E_\omega$ equals $\mathbb{C}$ while 
$E_\nu = (0)$ for any $\nu\in \Omega\setminus \{\omega\}$. 
There exists $y\in \Gamma(E)$ such that $y(\omega) =1$. 
Thus, 
$$\abs{y}(\nu) = \begin{cases}
0 &\qquad {\rm if\ }\nu\neq \omega\\
1 &\qquad {\rm if\ }\nu = \omega
\end{cases}$$ 
(which is clearly not continuous as $\omega$ is not an isolated point) and so
$\|y\|_e =0$ (in fact, $\|x\|_e =0$ for any $x\in X$). 
\end{eg}

\medskip

\begin{thm}
\label{ston}
Let $\Omega$ be a Stonian space and $X$ be a $C(\Omega)$-convex Banach module. 
Then $X_\re = X_\e$. 
Consequently, if $X$ comes from an (F)-Banach bundle, then $X$ is regular.
\end{thm}
\begin{prf}
By the argument of ``(i)$\Rightarrow$(ii)'' in Proposition \ref{equiv-loc-reg}, we see that for any $x\in X$ and 
$\epsilon > 0$, there exists $f\in C(\Omega)_+$ with $f\leq \abs{x}$ and $\|f\| = \|x\|_e - \epsilon$. 
Moreover, by the argument of ``(ii)$\Rightarrow$(iii)'' in Proposition \ref{equiv-loc-reg}, the map 
$\varphi: X(x) \rightarrow C(\Omega)$ given by $\varphi(g\cdot x) = gf$ is well defined. 
As $f \leq \abs{x}$ and $f$ is continuous, 
$$\|gf\| \ = \ \|gf\|_e \ \leq \ \|g\cdot x\|_e$$ 
for any $g\in C(\Omega)$. 
Therefore, $\varphi$ is a contraction from the semi-normed space $(X(x), \|\cdot \|_e)$ to $C(\Omega)$ and so it defines a contraction in $(X(x)_\e)^s$, also denoted by $\varphi$, such that $\|x_0\|_e \leq \|\varphi(x_0)\| + \epsilon$ (where $x_0$ is the image of $x$ in $X(x)_\e$). 
It is not hard to check that $X(x)_\e = X_\e(x_0)$. 
Thus, as $\Omega$ is Stonian, $\varphi$ extends to $\psi\in (X_\e)^s$ such that $\|\psi\| = \|\varphi\| \leq 1$ (by \cite[3.10]{KR}). 
Hence $X_\e$ is regular. 

On the other hand, suppose that $Y$ is a regular Banach $C(\Omega)$-bimodule and $\Phi\in \mathcal{L}_{C(\Omega)}(X,Y)$. 
Let $E$ and $F$ be (H)-Banach bundles over $\Omega$ such that $X = \Gamma(E)$ and $Y = \Gamma(F)$
(note that as $Y$ is regular, it is $C(\Omega)$-convex and such $F$ exists). 
Then $\Phi$ induces a Banach bundle map $\Psi: E \rightarrow F$. 
By Proposition \ref{equiv-loc-reg}, $\|y\| = \|y\|_e$ for any $y\in Y$. 
Thus for any $z\in X$, 
\begin{eqnarray*}
\|\Phi(z)\| 
& = & \|\Psi\circ z\| 
\ = \ \|\Psi\circ z\|_e 
\ = \ \inf \left\{ \sup_{\omega\in \Delta} \|\Psi(z(\omega))\|: \Delta {\rm\ is\ dense\ in \ }\Omega\right\}\\
& \leq & \|\Phi\| \cdot \inf \left\{ \sup_{\omega\in \Delta} \|z(\omega)\|: \Delta {\rm\ is\ dense\ in \ }\Omega\right\}
\ = \ \|\Phi\| \cdot \|z\|_e.
\end{eqnarray*}
Consequently, $\Phi$ factors through an element in $\mathcal{L}_{C(\Omega)}(X_\e,Y)$ uniquely. 
Since $X_\e$ is regular, it is the regularization of $X$. 
The second statement comes from the fact that $\|x\| = \|x\|_e$ if $X$ comes from an (F)-Banach bundle. 
\end{prf}

\medskip

This theorem and Remark \ref{t2>norm} show that if $\Omega$ is a Stonian space and 
$E$ is a Hausdorff (H)-Banach bundle over $\Omega$ with $X = \Gamma(E)$, 
then $\kappa_X$ is injective. 

\medskip

\begin{rem}
\label{ext}
Suppose that $\Omega$ is a compact Hausdorff space, $Y$ is a normed $C(\Omega)$-module and\ 
$n:Y \rightarrow C(\Omega)_+$\ satisfy the three conditions in \cite[p.47]{DG}. 
Then $n$ extends to the completion $\ti Y$ of $Y$ which also satisfies the same 
three conditions. 
In fact, $-n(y - z) \leq n(y) - n(z) \leq n(y - z)$ implies that 
$\|n(y)-n(z)\| \leq \|n(y - z)\| = \| y - z\|$ 
($y,z\in Y$). 
It is not hard to see that this gives a well-defined map from $\ti Y$ to $C(\Omega)_+$ (as $C(\Omega)_+$ 
is complete) which satisfies the three conditions in \cite[p.47]{DG}. 
\end{rem}

\medskip

\begin{thm}
\label{hy-ston}
Let $\Omega$ be a hyper-Stonian space and $X$ be a $C(\Omega)$-convex Banach module. 
Then $X_\e$ is the space of continuous sections of an (F)-Banach bundle. 
\end{thm}
\begin{prf}
Let $\{\mu_i\}_{i\in I}$ be a maximal family of positive normal measures on $\Omega$ with disjoint supports and 
let $\Xi_i$ be the support of $\mu_i$. 
Then $\Xi := \bigcup_{i\in I} \Xi_i$ is an open dense subset of $\Omega$ and $\{\mu_i\}_{i\in I}$ defines a 
Radon measure $\mu$ on $\Xi$ such that 
$C(\Omega) \cong L^\infty(\Xi, \mu)$ (see the argument of ``(i)$\Rightarrow$(ii)'' in \cite[III.1.18]{Tak}). 
Denote by ${\rm USC}(\Omega)_+$ the set of all upper semi-continuous functions from $\Omega$ to $\mathbb{R}_+$.
For any $h\in {\rm USC}(\Omega)_+$, we let $\psi(h)$ be the 
equivalence class of $h\!\mid\!_\Xi$ in $L^\infty(\Xi, \mu)$. 
We 
first show that 
$$\|h\|_e = \|\psi(h)\|_\infty.$$ 
Let $\Lambda\subseteq \Xi$ be a measurable set such that $\mu(\Lambda) =0$. 
Then $\mu_i(\Lambda) = 0$ for all $i\in I$ and so $\Lambda$ is nowhere dense in $\Omega$ (see \cite[III.1.15]{Tak}). 
The set $\Delta = \overline{\Lambda}\cup (\Omega\setminus \Xi)$ is closed and nowhere dense in $\Omega$.
Since $\Omega \setminus \Delta \subseteq \Xi\setminus \Lambda$, we see that 
$$\sup_{\omega \in \Omega \setminus \Delta} h(\omega)\ \ \leq\ \sup_{\omega\in \Xi \setminus \Lambda} h(\omega).$$ 
As $\Lambda$ is an arbitrary measure zero set, $\|h\|_e \leq \|\psi(h)\|_\infty$. 
Conversely, suppose that $\Delta$ is a closed nowhere dense subset of $\Omega$ and let 
$\Lambda = \Delta\cap \Xi$. 
Let $C$ be a compact subset of $\Lambda$. 
Then $C$ is also nowhere dense in $\Omega$ and $\mu_i(C) =0$ for any $i\in I$ (because of \cite[III.1.15]{Tak}). 
Therefore, $\mu(C) = \sum_{i\in I} \mu_i(C) = 0$ and so by the regularity of $\mu$, we have $\mu(\Lambda) = 0$. 
Since 
$$\sup_{\omega\in \Xi\setminus \Lambda} h(\omega)\ \ \leq \ \sup_{\omega\in \Omega\setminus \Delta} h(\omega)$$ 
and $\Delta$ is an arbitrary closed nowhere dense subset of $\Omega$, we see that 
$\|\psi(h)\|_\infty \leq \|h\|_e$. 
Now, as $g\mapsto \psi(g)$ is the canonical isomorphism from $C(\Omega)$ to $L^\infty(\Xi, \mu)$ (see the argument of \cite[III.1.18]{Tak}), the above shows that for any $h\in {\rm USC}(\Omega)_+$, there exists a unique $g\in C(\Omega)_+$ such that $g=h$ on an open and dense subset of $\Xi$ (and hence of $\Omega$). 
This induces a map 
$$\phi: {\rm USC}(\Omega)_+ \rightarrow C(\Omega)_+$$ 
such that $\|\phi(h)\| = \|h\|_e$.
For any $g\in C(\Omega)_+$ and $h\in USC(\Omega)_+$, it is not hard to see that $\phi(gh) = g\phi(h)$. 
If we set $n(x) := \phi(\abs{x})\in C(\Omega)_+$ ($x\in X$), then $\|x\|_e = \|n(x)\|$ 
(because $\abs{x}$ is upper semi-continuous). 
Hence, if $x,y\in X$ such that $\|x-y\|_e = 0$, then $\|(\abs{x-y})\!\mid\!_\Xi\|_\infty =0$ and so $\abs{x-y} = 0$ \emph{a.e.} on $\Xi$ which implies that $\abs{x} = \abs{y}$ on an open dense subset of $\Xi$. 
This shows that if $X_\e^0$ is the image of $X$ in $X_\e$, then $n$ induces 
$$\ti n: X_\e^0 \rightarrow C(\Omega)_+$$ 
which satisfies conditions (1)-(3) in \cite[p.47]{DG}. 
In fact, for any $g\in C(\Omega)$ and any $y\in X$, 
$$\ti n(g\cdot y_0) \ = \ \ti n((g\cdot y)_0)\ =\ \phi(\abs{g\cdot y})\ =\ \abs{g}\phi(\abs{y}) \ = \ 
\abs{g} \ti n(y_0)$$ 
(where $(g\cdot y)_0$ and $y_0$ are the images of $g\cdot y$ and $y$ respectively in $X_\e$). 
By Remark \ref{ext}, we see that $\ti n$ can be extended to $X_\e$. 
Thus, $X_\e$ is the space of continuous sections of an (F)-Banach bundle (see e.g. \cite[p.47-48]{DG}). 
\end{prf}

\medskip

\begin{rem}
We would like to thank the referee for informing us about \cite{Mag-dual} and for telling us that one can use the results in \cite[Section 6]{Mag-dual} to obtain the above theorem in an easier way. 
We left it to the readers to check the details. 
We decided to keep the proof as above because it is more elementary and our approach is completely different from the results in \cite{Mag-dual}.
\end{rem}

\medskip

\begin{defn}
Let $E$ and $E_c$ be respectively an (H)-Banach bundle and an (F)-Banach bundle over $\Omega$ and 
let $\Psi$ be a Banach bundle map from $E$ to $E_c$. 
Then $(E_c, \Psi)$ is called the \emph{continuous envelop of $E$} if any Banach bundle map from $E$ to any 
(F)-Banach bundle over $\Omega$ factors through $\Psi$ uniquely. 
\end{defn}

\medskip

It is not known if $(E_c, \Psi)$ always exists but it is 
the case when $\Omega$ is a hyper-Stonian space because of Theorems \ref{ston} and \ref{hy-ston}. 

\medskip

\begin{cor}
\label{ess=cont}
Suppose that $\Omega$ is a hyper-Stonian space, $E$ is an (H)-Banach bundle over $\Omega$ and $X = \Gamma(E)$. 
Then $E_c$ exists and $X_\e = \Gamma(E_c)$. 
\end{cor}

\medskip

\begin{cor}
\label{ess=reg=(F)}
Let $\Omega$ be a hyper-Stonian space and $X$ be a $C(\Omega)$-convex 
Banach module. 
The following are equivalent:

\smnoind
(i) $X = X_\e$. 

\smnoind
(ii) $X$ is regular. 

\smnoind
(iii) $X = \Gamma(E)$ for an (F)-Banach bundle over $\Omega$. 
\end{cor}

\medskip

\begin{rem}
As seen in Example \ref{eg-h}, not every (H)-Banach bundle over a hyper-Stonian space is an (F)-Banach bundle.
Therefore, regularity and pseudo-regularity do not coincide in this case. 
\end{rem}

\medskip
\noindent
School of Mathematical Sciences and LPMC,
Nankai University, 
Tianjin 300071,
China

\smnoind 
E-mail address: ckng@nankai.edu.cn

\end{document}